\documentclass[12pt]{article}
\usepackage{amsmath}
\usepackage{amssymb}
\usepackage{amsthm}
\usepackage{mathrsfs}
\usepackage{comment}
\usepackage[capitalise]{cleveref}
\crefformat{equation}{(#2#1#3)}
\usepackage{enumerate}
\usepackage{graphicx}
\usepackage[margin=1.5in]{geometry}

\newtheorem{thm}{Theorem}
\newtheorem*{thm*}{Theorem}
\newtheorem*{oq*}{Open Question 3.10}

\newtheorem*{cor*}{Corollary}
\newtheorem{conj}[thm]{Conjecture}

\theoremstyle{definition}

\newtheorem{qst}[thm]{Question}
\newtheorem*{qst*}{Question}
\newtheorem*{ua}{Ultrapower Axiom}

\makeatletter\@addtoreset{case}{thm}\makeatother
\makeatletter\@addtoreset{clm}{thm}\makeatother
\makeatletter\@addtoreset{sclm}{thm}\makeatother
\begin{document}
\title{Strong Compactness and the Ultrapower Axiom}
\author{Gabriel Goldberg}
\maketitle
\section{Some consequences of UA}

We announce some recent results relevant to the inner model problem. These results involve the development of an abstract comparison theory from a hypothesis called the Ultrapower Axiom, a natural generalization to larger cardinals of the assumption that the Mitchell order on normal measures is linear. Roughly, it says that any pair of ultrapowers can be ultrapowered to a common ultrapower. The formal statement is not very different. Let \(\text{Uf}\) denote the class of countably complete ultrafilters.

\begin{ua}
For any \(U_0,U_1\in \textnormal{Uf}\), there exist \(W_0\in \textnormal{Uf}^{M_{U_0}}\) and \(W_1\in \textnormal{Uf}^{M_{U_1}}\) such that \[M^{M_{U_0}}_{W_0} = M^{M_{U_1}}_{W_1}\] and \[j^{M_{U_0}}_{W_0}\circ j_{U_0} = j^{M_{U_1}}_{W_1}\circ j_{U_1}\]
\end{ua}

The Ultrapower Axiom, which we will refer to from now on as UA, holds in all known inner models. This is a basic consequence of the methodology used to build these models, the comparison process. For example, assuming there is a proper class of strong cardinals and \(V = \text{HOD}\), UA follows from Woodin's principle Weak Comparison \cite{Woodin}, which is an immediate consequence of the basic comparison lemma used to construct and analyze canonical inner models. These hypotheses are just an artifact of the statement of Weak Comparison, and in fact the methodology of inner model theory simply cannot produce models in which UA fails (although of course it can produce models in which \(V\neq \text{HOD}\) and there are no strong cardinals). The following is therefore a precise test question for the inner model program.

\begin{qst}\label{con}
Is UA consistent with a supercompact cardinal?
\end{qst}

Conditioned on iterability, there are canonical inner models reaching the finite levels of supercompactness, and these models satisfy UA. It would be strange, then, if UA was refuted by the large cardinal hypothesis asserting the existence of a cardinal \(\kappa\) that is \(\kappa^{+\omega}\)-supercompact, and perhaps even stranger if UA were consistent with this hypothesis but ruled out by some other large cardinal hypothesis below a supercompact. We therefore think that the answer to \cref{con} is probably yes. Results of Woodin strongly suggest that if UA is consistent with a supercompact cardinal then it is consistent with all large cardinals. More formally, if there is a weak extender model for supercompactness that satisfies UA, then this model absorbs all large cardinals. This is our justification for stating \cref{ghod}, whose hypotheses might otherwise look quite speculative.

UA turns out to answer many questions in the structure theory of countably complete ultrafilters. More surprisingly, in the presence of a single strongly compact cardinal, UA begins to answer many set theoretic questions that on the face of it have nothing to do with ultrafilters.

\begin{thm}[UA]
\textnormal{GCH} holds above the least strongly compact cardinal.
\end{thm}

That is, if \(\kappa\) is the least strongly compact cardinal, then for all \(\lambda\geq \kappa\), \(2^\lambda = \lambda^+\). For singular strong limit cardinals \(\lambda\), this is due to Solovay and provable in ZFC.

Another class of results show that \(\textnormal{HOD}\) is quite close to \(V\) assuming UA.

\begin{thm}[UA]\label{hod}
If there is a strongly compact cardinal, then \(V\) is a generic extension of \(\textnormal{HOD}\).
\end{thm}

One cannot prove \(V = \text{HOD}\) outright from UA since UA is preserved under generic extensions below the least measurable. But one can show the forcing is not much larger than the least strongly compact \(\kappa\): the proof shows that \(V = \textnormal{HOD}_X\) for any \(X\subseteq \kappa\) such that \(V_\kappa\subseteq \text{HOD}_X\).

\begin{thm}[UA]
If there is a supercompact cardinal \(\kappa\), then \(\textnormal{HOD}\) is a weak extender model for \(\kappa\) is supercompact.
\end{thm}

The proof of \cref{hod} has the following corollary. Here the {\it mantle}, denoted \(\mathbb M\), is the intersection of all forcing grounds, and the {\it generic} HOD, denoted \(\textnormal{gHOD}\), is the intersection of \(\text{HOD}^{N}\) over all generic extensions \(N\). Both structures are inner models of ZFC (the proof that the mantle satisfies the Axiom of Choice uses the Downward Directed Grounds Hypothesis, a theorem of Usuba).

\begin{thm}[UA]
Suppose there is a proper class of supercompact cardinals. Then \(\mathbb M = \textnormal{gHOD}\). 
\end{thm}

Using Usuba's theorem, we also have the following theorem (whose hypothesis is of course well past \cref{con}).

\begin{thm}[UA]\label{ghod}
Suppose there is a hyperhuge cardinal. Then the mantle satisfies \(V = \textnormal{gHOD}\).
\end{thm}

We remark that this is (trivially) stronger than saying that the mantle satisfies \(V = \text{HOD}\). 

Probably the most interesting consequences of UA exhibit a very natural relationship between comparison and supercompactness. This suggests to us that the answer to \cref{con} is yes.

\begin{thm}[UA]\label{strongsuper}
The least strongly compact cardinal is supercompact.
\end{thm}

One can think of this as proving the ``canonical equiconsistency" of strongly compact and supercompact cardinals: if there is a canonical model with a strongly compact cardinal, there is a canonical model with a supercompact cardinal.

An obvious question is left open by the proof of \cref{strongsuper}:

\begin{qst*}[UA]
Is the second strongly compact cardinal supercompact?
\end{qst*}

We mention a local version of \cref{strongsuper} assuming GCH in the region just beyond where the inner models have been (conditionally) constructed, which yields some insight into the inner model problem at this level.

\begin{thm}[UA + GCH]\label{antiinner}
The least cardinal \(\kappa\) that is \(\kappa^{+\omega}\)-strongly compact is \(\kappa^{+\omega}\)-supercompact.
\end{thm}

The impetus for proving \cref{antiinner} was an argument due to Woodin that in a nonstrategic extender model, the least \(\kappa^{+\omega}\)-strongly compact is {\it not} \(\kappa^{+\omega}\)-supercompact. Therefore along with Woodin's work, \cref{antiinner} seems to show that nonstrategic extender models do not reach the infinite levels of strong compactness at all, presumably because of a failure of iterability. We can say with more certainty that Woodin's One Generator Framework cannot produce nonstrategic extender models at the infinite levels of strong compactness. This came as a bit of a surprise since the upper bound identified by Woodin \cite{Woodin} suggested
the nonstrategic extender hierarchy could reach the level of supercompact (but not the level of a supercompact with a measurable above).

\section{Conclusion}
There is enormous structure above the least strongly compact cardinal assuming UA. The cluster of theorems around this point fits perfectly with Woodin's isolation of the least supercompact cardinal as a critical juncture for inner model theory. The emerging theory relating strong compactness and supercompactness under UA suggests the consistency of UA with supercompact cardinals and beyond, and therefore gives some evidence that the inner model program will succeed. 

Assuming this is the correct picture, it is natural to wonder whether UA itself might axiomatize the ultimate inner model above the least supercompact cardinal \(\kappa\). Maybe UA answers everything ``above \(\kappa\)." We make this a precise conjecture. UA can be preserved by forcing in gaps between countably complete uniform ultrafilters, which exist only below the least supercompact cardinal. The idea is that perhaps this is the only obstacle preventing UA from settling all questions. 

\begin{conj}\label{recovery}
Assume the Ultrapower Axiom, the Ground Axiom, and the existence of a supercompact cardinal. Then \(V = \textnormal{Ultimate } L\).
\end{conj}

\bibliography{seedorder}{}
\bibliographystyle{unsrt}
\end{document}